\documentclass[11pt, draft]{article}
\usepackage{amsfonts}

\newtheorem{theorem}{\bf Theorem}[section]
\newtheorem{lemma}{\bf Lemma}[section]
\newenvironment{proof}{{ \medskip \noindent \bf Proof: }}{ \hfill
                              \rule{2mm}{3mm} \newline \medskip}

\title{The Dimension of Quasi-Homogeneous Planar Linear Systems With Multiplicity Four}
\author{James Seibert}

\begin{document}

\maketitle

\begin{abstract}  A linear system of plane curves satisfying multiplicity conditions at points in general position is called special if the dimension is larger than the expected dimension.  A (-1) curve is an irreducible curve with self intersection -1 and genus zero.  The Harbourne-Hirschowitz Conjecture is that a linear system is special only if a multiple of some fixed (-1) curve is contained in every curve of the linear system.  This conjecture is proven for linear systems with multiplicity four at all but one of the points.  \end{abstract}

\section{Introduction}
Consider $n+1$ points $p_0,p_1, \dots p_n$ in the projective plane ${\mathbb P}$.  Let $\mathcal{L}$ be the linear system of plane curves of degree $d$ with multiplicity $m_i$ at point $p_i$ for $0\leq i\leq n$.  When all the multiplicities $m_i$ for $i\geq1$ are equal to some value $m$ the system is called {\it quasi-homogeneous} and will be written as $\mathcal{L}=\mathcal{L}(d, m_0, n, m)$.  The space of curves of degree $d$ has dimension $\frac{d(d+3)}{2}$, and a point of multiplicity $k$ imposes $\frac{k(k+1)}{2}$ linear conditions.  This leads us to define the {\it virtual dimension} of $\mathcal{L}$ to be
$$v= \frac{d(d+3)}{2} - \frac{m_0(m_0+1)}{2} - n \frac{m(m+1)}{2}.$$
Of course the dimension cannot be less than -1, so we define the {\it expected dimension}
$$e= \max \{ -1, v \}.$$
There is a dense open subset in the parameter space of $n+1$-tuples of points on which the actual dimension of $\mathcal{L}$ achieves a minimum.  This we call the {\it dimension} of $\mathcal{L}$, and we write $\dim \mathcal{L}=\ell$.
We know that $\ell \geq e$.  When $\ell > e$ the conditions imposed by the points are dependent and the system is called {\it special}.  When $\ell = e$ the system is {\it non-special}.

One of the methods used to find the dimension of a particular linear system is to reduce the data, when possible, by using standard quadratic transformations of the plane, also known as Cremona transformations.  These transformations do not change the actual dimension of a system, and this is useful when the transformed system has known dimension.  The words ``reduce'' and ``Cremona reduce'' refer to this process throughout this paper.

We define the {\it intersection number} (when $n^\prime \leq n$),
$$\mathcal{L}(d, m_0, n, m) \cdot \mathcal{L}(d^\prime, m_0^\prime, n^\prime, m^\prime) = d d^\prime - m_0 m_0^\prime - n^\prime m m^\prime.$$
The {\it self-intersection number} is then defined in the natural way,
$$\mathcal{L}^2 = d^2 -m_0^2 -nm^2.$$
Pl\"ucker's formula for the {\it genus} of curves in $\mathcal{L}$ is
$$g_{\mathcal{L}} = \frac{(d-1)(d-2)}{2} - \frac{m_0(m_0-1)}{2} - n\frac{m(m-1)}{2}.$$
A {\it (-1) curve} is an irreducible curve in a system with self intersection -1 and with genus zero.  A multiple (-1) curve is contained in all known examples of special linear systems.

A (-1) curve contained in a quasi-homogeneous linear system $\mathcal{L}$ must either be quasi-homogeneous itself, or be quasi-quasi-homogeneous.  By quasi-quasi-homogeneous we refer to a curve $C_k$ with multiplicity $m_0$ at point $p_0$, multiplicity $m_1$ at one of the points $p_k$ for $1\leq k \leq n$, and multiplicity $m$ at the rest of the points $p_1$ through $p_n$ (except $p_k$).  If $\mathcal{L}$ contains one such curve $C_k$, then it must contain the sum of the curves $\sum_{i=1}^n C_k$.  This sum is quasi-homogeneous, and is referred to as a quasi-homogeneous (-1) configuration.

A linear system $\mathcal{L}$ is {\it (-1) special} if there are (-1) curves $A_1, \dots, A_r$ 
such that $\mathcal{L} \cdot A_j = - N_j$ with $N_j \geq 1$ for every $j$ and $N_j \geq 2$ for 
some $j$, with the residual system $\mathcal{M} = \mathcal{L} - \sum_j N_j  A_j$ having 
non-negative virtual dimension and non-negative intersection with every (-1) curve $A_j$.  The 
(-1) curves $A_j$ must be pairwise disjoint. 

It is known that every (-1) special system is special. 
The main conjecture is that every special system is (-1) special.  We begin by classifying the
(-1) special linear systems with $m=4$.

\section{(-1) Special Systems}
Suppose that $\mathcal{L} (d, m_0, n, 4)$ is a (-1) special linear system.  Then $\mathcal{L}$ 
must be one of two possible forms.  The first possibility is that $$\mathcal{L} = \mathcal{M} + 
N \cdot C,$$ where $v(\mathcal{M}) \geq 0$ and $\mathcal{M} \cdot C = 0$.  In this case $N$ is 
2, 3, or 4, and the (-1) curve or configuration $C = \mathcal{L}(\delta, \mu_0, n, 1)$, or $N$ 
is 2 and $C = \mathcal{L}(\delta, \mu_0, n, 2)$.  The second possibility is that $$\mathcal{L} 
= 2C_1 + 2C_2.$$  In this case $C_1$ and $C_2$ must be disjoint (-1) curves or configurations.
The proof of the following Lemma is not difficult.
\begin{lemma}
The following is a complete list of (-1) curves and configurations with $m\leq2$.

$$\vbox{

\tabskip = 15pt
\halign{
        \hfil#&
        \hfil#&
        #\hfil

\cr
        $\mathcal{L}(2, 0, 5, 1)$ & &  \cr
        $\mathcal{L}(e, e-1, 2e, 1)$ & $e\geq1$ &  \cr
        $\mathcal{L}(e, e, e, 1)$ & $e\geq1$ & compound for $e>2$ \cr
        $\mathcal{L}(6, 3, 7, 2)$ & &   \cr
        $\mathcal{L}(3, 0, 3, 2)$ & & compound  \cr }
   
}$$
\end{lemma}  
Note that the only disjoint curves on this list are $\mathcal{L}(e, e-1, 2e, 1)$ and 
$\mathcal{L}(2e, 2e, 2e, 1)$.   
We may now classify (-1) special curves with $m=4$.
\begin{theorem} \label{mosp} The quasi homogeneous (-1) special systems with $m=4$ are those 
systems $\mathcal{L}(d, m_0, n, 4)$ on the following list:
\end{theorem}
\begin{tabbing}
\hspace{.5in} \= $d-m_0$ \hspace{.2in} \= \hspace{1.8in} \= \\
   \>   0 \> $\mathcal{L}(d, d, e, 4)$ \> $d \geq 4e \geq 4$ \\ 
   \>   1 \> $\mathcal{L}(d, d-1, e, 4)$ \> $d \geq \frac{7}{2}e \geq 1$ \\
   \>   2 \> $\mathcal{L}(d, d-2, e, 4)$ \> $d\geq \frac{9e+1}{3} \geq \frac{10}{3}$ \\
   \>   2 \> $\mathcal{L}(6e, 6e-2, 2e, 4)$ \> $e\geq1$ \\
   \>   3 \> $\mathcal{L}(5e, 5e-3, 2e, 4)$ \> $e\geq1$ \\
   \>   3 \> $\mathcal{L}(5e+1, 5e-2, 2e, 4)$ \> $e\geq1$ \\
   \>   4 \> $\mathcal{L}(4e, 4e-4, 2e, 4)$ \> $e\geq1$ \\
   \>   4 \> $\mathcal{L}(4e+1, 4e-3, 2e, 4)$ \> $e\geq1$ \\
   \>   4 \> $\mathcal{L}(4e+2, 4e-2, 2e, 4)$ \> $e\geq1$ \\
   \>   5 \> $\mathcal{L}(5, 0, 2, 4)$ \\
   \>   5 \> $\mathcal{L}(6, 1, 2, 4)$ \\
   \>   5 \> $\mathcal{L}(8, 3, 4, 4)$ \\
   \>   5 \> $\mathcal{L}(9, 4, 4, 4)$ \\
   \>   5 \> $\mathcal{L}(12, 7, 6, 4)$ \\
   \>   5 \> $\mathcal{L}(15, 10, 8, 4)$ \\
   \>   6 \> $\mathcal{L}(6, 0, 2, 4)$ \\
   \>   6 \> $\mathcal{L}(6, 0, 3, 4)$ \\
   \>   6 \> $\mathcal{L}(8, 2, 4, 4)$ \\
   \>   6 \> $\mathcal{L}(12, 6, 7, 4)$ \\
   \>   7 \> $\mathcal{L}(9, 2, 5, 4)$ \\
   \>   8 \> $\mathcal{L}(8, 0, 5, 4)$ \\
   \>   8 \> $\mathcal{L}(9, 1, 5, 4)$ \\
   \>   9 \> $\mathcal{L}(9, 0, 5, 4)$ \\
 \end{tabbing}
\begin{proof}
Suppose that $\mathcal{L} = \mathcal{M} + NC$ with $N = 2, 3$, or 4.  We begin with $N=4$.

Suppose that $\mathcal{L}(2, 0, 5, 1)$ splits off 4 times from $\mathcal{L}(d, m_0, 5, 4)$.  
Then the residual system $\mathcal{M} = \mathcal{L}(d-8, m_0, 5, 0)$.  $\mathcal{M}$ is 
disjoint from $\mathcal{L}(2, 0, 5, 1)$ only if $d=8$, and so $v(\mathcal{M}) = 0$ implies 
that $m_0 = 0$.  Thus  we have the class $4 \cdot \mathcal{L}(2, 0, 5, 1) = 
\mathcal{L}(8, 0, 5, 4)$ with virtual dimension $-6$, but dimension $0$.

Suppose that $A = \mathcal{L}(e, e-1, 2e, 1)$ splits off 4 times from 
$\mathcal{L}(d, m_0, 2e, 4)$.  Then $\mathcal{M} = \mathcal{L}(d-4e, m_0-4e+4, 2e, 0)$.  
We require that 
$$\mathcal{M} \cdot A = e(d-m_0-8) +m_0+4 = 0$$
so that $m_0 > d-8$.  Then $v(\mathcal{M}) \geq 0$ implies that $m_0 \leq d-4$.  There are four 
possibilities for $m_0$.

$m_0 = d-7$:  $\mathcal{M} \cdot A = e(d-m_0-8) +m_0+4 = 0$ implies that $e = d-3$.  Then 
$v(\mathcal{M}) = -12d +39 \geq 0$ implies that $d \leq 3$.  This, however,  makes $m_0 < 0.$ 

$m_0 = d-6$:  $\mathcal{M} \cdot A = e(d-m_0-8) +m_0+4 = 0$ implies that $d = 2e + 2$.  But 
$d \geq 4e$ and $e\geq1$, so $e=1$ and $d=4$.  Again this makes $m_0 <0$.

$m_0 = d-5$:  $\mathcal{M} \cdot A = e(d-m_0-8) +m_0+4 = 0$ implies that $d = 3e+1$.  But 
$d \geq 4e$ and $e\geq1$, so $e=1$ and $d=4$.  This makes $m_0 = -1$.

$m_0 = d-4$:  $\mathcal{M} \cdot A = e(d-m_0-8) +m_0+4 = 0$ implies that $d = 4e$.  This gives 
the system $\mathcal{L}(4e, 4e-4, 2e, 4)$, for all $e\geq1$.

Suppose that $A = \mathcal{L}(e, e, e, 1)$ splits off 4 times from 
$\mathcal{L}(d, m_0, e, 4)$.  Then $\mathcal{M} = \mathcal{L}(d-4e, m_0-4e, 2e, 0)$.  We 
require that $$\mathcal{M} \cdot A = e(d-4e) +e(m_0-4e) = 0,$$ so that $d=m_0$.  Now 
$v(\mathcal{M}) \geq 0$ whenever $d\geq 4e$.  This gives the system $\mathcal{L}(d, d, e, 4)$ 
whenever $d\geq4e \geq 4$.

This completes the $N=4$ analysis.

Suppose that $A = \mathcal{L}(2, 0, 5, 1)$ splits off 3 times from $\mathcal{L}(d, m_0, 5, 4)$.
Then $\mathcal{M} = \mathcal{L}(d-6, m_0, 5, 1)$, and $\mathcal{M} \cdot A = 2d -17$.  This can 
never be zero, so this case cannot occur.

Suppose that $A = \mathcal{L}(e, e-1, 2e, 1)$ splits off 3 times from 
$\mathcal{L}(d, m_0, 2e, 4)$.  Then $\mathcal{M} = \mathcal{L}(d-3e, m_0-3e+3, 2e, 1)$.  We 
require that 
$$\mathcal{M} \cdot A = e(d-m_0-8)+m_0 +3 = 0.$$
This forces $m_0 > d-8$.  Then $v(\mathcal{M}) \geq 0$ requires $m_0 \leq d-3$, and there are 
five possibilities for $m_0$.

$m_0 = d-7$:  $\mathcal{M} \cdot A = d-e-4 = 0$ implies that $m_0 = e-3$, but $m_0 \geq 3e-3$.  
This case cannot occur.

$m_0 = d-6$:  $\mathcal{M} \cdot A = d-2e-3 = 0$ implies that $m_0 = 2e-3$, but 
$m_0 \geq 3e-3$.  This case cannot occur.

$m_0 = d-5$:  $\mathcal{M} \cdot A = d-3e-2 = 0$ implies that $d = 3e+2$ and $m_0 = 3e-3$.  
In this case $\mathcal{M} = \mathcal{L}(2, 0, 2e, 1)$ and $v(\mathcal{M}) \geq 0$ forces $e$ to 
be 1 or 2.  These are the systems $\mathcal{L}(5, 0, 2, 4)$ and $\mathcal{L}(8, 3, 4, 4)$.

$m_0 = d-4$:  $\mathcal{M} \cdot A = d-4e-1 = 0$, so $d=4e+1$ and $m_0 = 4e -3$.  Now 
$\mathcal{M} = \mathcal{L}(e+1, e, 2e, 1)$ has virtual dimension 2 for all $e$.  Therefore we 
have the (-1) special system $\mathcal{L}(4e+1, 4e-3, 2e, 4)$ for all $e\geq1$.

$m_0 = d-3$:  $\mathcal{M} \cdot A = d-5e = 0$, so $d=5e$ and $m_0 = 5e -3$.  Now 
$\mathcal{M} = \mathcal{L}(2e, 2e, 2e, 1)$ has virtual dimension 0 for all $e$.  This gives the 
system $\mathcal{L}(5e, 5e-3, 2e, 4)$ for all $e\geq1$.

Suppose that $A = \mathcal{L}(e, e, e, 1)$ splits off 3 times from $\mathcal{L}(d, m_0, e, 4)$. 
Then $\mathcal{M} = \mathcal{L}(d-3e, m_0-3e, 2e, 1)$.  We require that $$\mathcal{M} \cdot A 
= e(d-m_0 -1) = 0,$$ so that $m_0 = d-1$.  Now $v(\mathcal{M}) = 2d -7e \geq 0 $ if $d \geq 
\frac {7}{2}e$.  This gives the systems $\mathcal{L}(d, d-1, e, 4)$ for every 
$d\geq \frac {7}{2}e \geq \frac{7}{2}$.

This completes the $N=3$ analysis.

Suppose that $A = \mathcal{L}(2, 0, 5, 1)$ splits off twice from $\mathcal{L}(d, m_0, 5, 4)$.  
Then $\mathcal{M} = \mathcal{L}(d-4, m_0, 5, 2)$, and $\mathcal{M} \cdot A = 2d -18$, so $d=9$. 
$v(\mathcal{M}) \geq 0$ implies that $0\leq M-0 \leq 2$.  This gives the systems 
$\mathcal{L}(9, 0, 5, 4)$, $\mathcal{L}(9, 1, 5, 4)$, and $\mathcal{L}(9, 2, 5, 4)$.

Suppose that $A = \mathcal{L}(e, e-1, 2e, 1)$ splits off twice from 
$\mathcal{L}(d, m_0, 2e, 4)$.  Then $\mathcal{M} = \mathcal{L}(d-2e, m_0-2e+2, 2e, 2)$ and 
$\mathcal{M} \cdot A = e(d-m_0-8) + m_0 +2$.  This forces $m_0 > d-8$.  We also have 
$m_0 \leq d-2$ so that  $v(\mathcal{M}) \geq 0$.  There are six possibilities for $m_0$.

$m_0 = d-7$:  $\mathcal{M} \cdot A = d-e-5 = 0$, so $d=e+5$ and $m_0 = e -2$.  But 
$m_0 \geq 2e -2$, so this case cannot occur.

$m_0 = d-6$:  $\mathcal{M} \cdot A = d-2e-4 = 0$, so $d=2e+4$ and $m_0 = 2e -2$.  
Now $v(\mathcal{M}) = 14 - 6e \geq 0$ implies that $e$ is 1 or 2.  This leads to the systems 
$\mathcal{L}(6, 0, 2, 4)$ and $\mathcal{L}(8, 2, 4, 4)$.

$m_0 = d-5$:  $\mathcal{M} \cdot A = d-3e-3 = 0$, so $d=3e +3$ and $m_0 = 3e-2$.  
Then $v(\mathcal{M}) = 9 - 2e \geq 0$ implies that $e$ is 1, 2, 3, or 4.  These systems are 
$\mathcal{L}(6, 1, 2, 4)$, $\mathcal{L}(9, 4, 4, 4)$, $\mathcal{L}(12, 7, 6, 4)$, and  
$\mathcal{L}(15, 10, 8, 4)$.

$m_0 = d-4$:  $\mathcal{M} \cdot A = d-4e-2 = 0$, so $d=4e +2$ and $m_0 = 4e-2$.  In this case 
$v(\mathcal{M}) = 5$ for all $e$, giving the systems $\mathcal{L}(4e+2, 4e-2, 2e, 4)$ for all 
$e\geq1$.

$m_0 = d-3$:  $\mathcal{M} \cdot A = d-5e-1 = 0$, so $d=5e +1$ and $m_0 = 5e-2$.  This time 
$v(\mathcal{M}) = 2$ for all $e$, giving the systems $\mathcal{L}(5e+1, 5e-2, 2e, 4)$ for all 
$e\geq1$.

$m_0 = d-2$:  $\mathcal{M} \cdot A = d-6e = 0$, so $d=6e$ and $m_0 = 6e-2$, but 
$v(\mathcal{M}) = -2$ for all $e$.  This means that $\mathcal{L}(e, e-1, 2e, 1)$ cannot split 
off by itself in this case.  We will see later that this case does occur when 
$\mathcal{L}(e, e-1, 2e, 1)$ and $\mathcal{L}(2e, 2e, 2e, 1)$ both split.

Suppose that $A = \mathcal{L}(e, e, e, 1)$ splits off twice from $\mathcal{L}(d, m_0, e, 4)$.  
Then $\mathcal{M} = \mathcal{L}(d-2e, m_0-2e, e, 2)$ and $\mathcal{M} \cdot A = e(d-m_0-2)$.  
This forces $m_0 = d-2$.  Now $v(\mathcal{M}) = 3d -9e -1 \geq 0$ occurs whenever 
$d\geq \frac{9e+1}{3}$.  This gives a (-1) special system $\mathcal{L}(d, d-2, e, 4)$ whenever 
$d \geq \frac{9e+1}{3} \geq \frac{10}{3}$.

Suppose that $A = \mathcal{L}(6, 3, 7, 2)$ splits off twice from $\mathcal{L}(d, m_0, 7, 4)$.  
Then $\mathcal{M} \cdot A = 6d -3m_0-54$ and $m_0 = 2d - 18$.  Now $v(\mathcal{M}) \geq 0$ 
only if $d \leq 12$, but $d\geq 12$ since $A$ splits off twice.  Therefore $d=12$ and the 
system is $\mathcal{L}(12, 6, 7, 4) = 2A$

Suppose that $A = \mathcal{L}(3, 0, 3, 2)$ splits off twice from $\mathcal{L}(d, m_0, 3, 4)$.  
Then   $\mathcal{M} \cdot A = 3(d-6) = 0$ implies that $d=6$, and $v(\mathcal{M}) \geq 0$ 
implies that $m_0 = 0$, giving the system $\mathcal{L}(6, 0, 3, 4) = 2A$.

This completes the analysis of the case where one curve splits off twice.  The final 
possibility is that $\mathcal{L} = \mathcal{M} + 2N_1 +2N_2$, where $N_1$ and $N_2$ are 
disjoint (-1) curves or configurations.  As mentioned previously, $\mathcal{L}(e, e-1, 2e, 1)$ 
and $\mathcal{L}(2e, 2e, 2e, 1)$ are the only possibilities for $N_1$ and $N_2$.  In this 
situation, the residual system is $\mathcal{M} = \mathcal{L}(d-6e, m_0-6e+2, 2e, 0)$.  $N_1$ 
and $N_2$ intersect $\mathcal{M}$ as 
$$\mathcal{M} \cdot \mathcal{L}(2e, 2e, 2e, 1) = 2e(d-m_0) -4e = 0.$$
$$\mathcal{M} \cdot \mathcal{L}(e, e-1, 2e, 1) = e(d-m_0-8)+m_0+2 = 0.$$
From the first equation we get that $m_0 = d-2$.  From the second we then have $d = 6e$.  
Now $v(\mathcal{M}) \geq 0$ implies that $m_0 = 6e -2$.  This leads to 
the (-1) special systems $2N_1 + 2N_2 = \mathcal{L}(6e, 6e-2, 2e, 4)$ for 
every $e\geq1$.  

\end{proof}

\section{Large $m_0$} \label{largem_0}

The cases with $m_0 \geq d-5$ are dealt with by the following lemmas from 
\cite{miranda} (modified to the case $m=4$).

\begin{lemma} Let $\mathcal{L} = \mathcal{L}(d, d-4, n, 4)$, with $d\geq4$.  
Write $d=4q+\mu$ with $0\leq \mu \leq 3$, and $n=2h + \epsilon$, with 
$\epsilon \in \{0, 1 \}$.  Then the system $\mathcal{L}$ is special if and 
only if $q=h$, $\epsilon = 0$ and $\mu \leq 2$.
\end{lemma}

In particular, $\mathcal{L}(d, d-4, n, 4)$ is special if and only if it is 
one of the following types:
$$\vbox{
\tabskip = 15pt
\halign{
        \hfil#&
        #\hfil
\cr
        $\mu = 0$ & $\mathcal{L}(4q, 4q-4, 2q, 4)$  \cr
	$\mu = 1$ & $\mathcal{L}(4q+1, 4q-3, 2q, 4)$ \cr
	$\mu = 1$ & $\mathcal{L}(4q+2, 4q-2, 2q, 4)$ \cr }
}$$
These agree with the list in theorem \ref{mosp}.

\begin{lemma} Let $\mathcal{L} = \mathcal{L}(d, d-4 +k, n, 4)$ with 
$k\geq1$, and let 
$$\mathcal{L}^\prime = \mathcal{L}(d-kn, d-kn-4+k, n, 4-k).$$
Then $\dim \mathcal{L} = \dim \mathcal{L}^\prime$ and $\mathcal{L}$ is 
non-special unless either

(a) $k\geq2$ and $\mathcal{L}^\prime$ is nonempty and non-special, or \\ 
\indent (b) $\mathcal{L}^\prime$ is special.
\end{lemma}

Again, we compare these results with the list in Theorem \ref{mosp}.  
If $k=1$, $\mathcal{L} = \mathcal{L}(d, d-3, n, 4)$ and 
$\mathcal{L}^\prime = \mathcal{L}(d-n, d-n-3, n, 3)$, a system with 
$m=3$.  $\mathcal{L}$ is special if and only if $\mathcal{L}^\prime$ is 
special.  The special quasi homogeneous systems with $m\leq3$ are 
classified in \cite{miranda}.  $\mathcal{L}^\prime$ is special if and only 
if it is of the form $\mathcal{L}(3e, 3e-3, 2e, 3)$, or 
$\mathcal{L}(3e +1, 3e-2, 2e, 3)$.  These lead to systems 
$\mathcal{L}(5e, 5e-3, 2e, 4)$ and $\mathcal{L}(5e+1, 5e-2, 2e, 4)$.  These 
appear on the list in Theorem \ref{mosp}, and they are the only classes on 
the list with $d-m_0 = 3$.

If $k=2$, $\mathcal{L} = \mathcal{L}(d, d-2, n, 4)$ and 
$\mathcal{L}^\prime = \mathcal{L}(d-2n, d-2n-2, n, 2)$.  
$\mathcal{L}^\prime$ is nonempty and non-special if $d \geq \frac{9n+1}{3}$.
  $\mathcal{L}^\prime$ is special if and only if it has the form 
$\mathcal{L}(2e, 2e-2, 2e, 2)$.  This gives 
$\mathcal{L} = \mathcal{L}(6e, 6e-2, 2e, 4)$.  
These agree with the list in Theorem \ref{mosp}.

If $k=3$, $\mathcal{L} = \mathcal{L}(d, d-1, n, 4)$ and 
$\mathcal{L}^\prime = \mathcal{L}(d-3n, d-3n-1, n, 1)$.  
$\mathcal{L}^\prime$ is non-special since $m=1$ and nonempty if 
$d \geq \frac{7n}{2}$.  If $k=4$, $\mathcal{L} = \mathcal{L}(d, d, n, 4)$ 
and $\mathcal{L}^\prime = \mathcal{L}(d-4n, d-4n, n, 0)$.  
$\mathcal{L}^\prime$ is non-special since $m=0$ and nonempty if 
$d\geq 4n$.  These conditions match those on the list in Theorem 
\ref{mosp}.

\begin{lemma} Let $\mathcal{L} = \mathcal{L}(d, d-5, n, 4)$ with $d\geq5$.  
Write $d=3q+\mu$ with $0 \leq \mu \leq 2$, and $n=2h+\epsilon$ with 
$\epsilon \in \{0, 1\}$.  Then the system $\mathcal{L}$ is non-special 
unless

(a) $q=h+1$, $\mu = \epsilon = 0$ and $h \leq 4$, or \\ \indent (b) $q=h$, 
$\epsilon = 0$, and $4q \leq \mu(\mu+3)$. \end{lemma}

Case (a) gives the systems $\mathcal{L}(3h+3, 3h-2, 2h, 4)$ where 
$1 \leq h \leq 4$.  
These appear on the list in Theorem \ref{mosp}.  In case (b), $\mu = 0$ 
and $\mu = 1$ both force $d\leq5$, so we are left with $\mu = 2$ and 
$q = 1$ or $2$.  These are the systems $\mathcal{L}(5, 0, 2, 4)$ and 
$\mathcal{L}(8, 3, 4, 4)$.  These account for all the systems with 
$d-m_0 = 5$ on the list in Theorem \ref{mosp}.

The main conjecture holds for quasi homogeneous systems with $m=4$ and 
$m_0 \geq d-5$.

\section{The Degeneration}

To finish the classification of special systems we use 
the degeneration of the plane described in detail by Ciliberto and Miranda 
\cite{miranda}.  The general plan relies on the fact that if the points are in special position,
the dimension of the system can only increase.  We attempt to find a special position for the
points that allows us to calculate the dimension, but such that the dimension is not greater 
than the expected dimension.  We are able to do so because the degeneration affords us a 
great deal of flexibility.

Briefly, we consider $V ={\mathbb P}^2 \times {\mathbb A}^1$ and $X$, the
blow-up of this three-fold along a line $L$ in ${\mathbb P}^2 \times \{ 0 \}$.  We have the projections
$p_1 : V \rightarrow {\mathbb A}^1$ and $p_2 : V \rightarrow {\mathbb P}^2$, the blowup map
$f : X \rightarrow V$, and the compositions $\pi_1 = p_1 \circ f : X \rightarrow 
{\mathbb A}^1$ and $\pi_2 = p_2 \circ f : X \rightarrow {\mathbb P}^2$.  Let $X_t$ be the 
fiber of $\pi_1$ over $t$ in ${\mathbb A}^1$.  Then $X_t \cong {\mathbb P}^2$ if 
$t \neq 0$.
In $X_0$, the degeneration produces two surfaces, a plane 
${\mathbb P} = {\mathbb P}^2$ and a Hirzebruch surface 
${\mathbb F} = {\mathbb F}_1$, joined transversely along a curve 
$R$.  $R$ is the line $L$ in ${\mathbb P}$ and the exceptional divisor $E$ 
of ${\mathbb F}$.  

The Picard group of $X_0$ is the fibered product of Pic(${\mathbb P}$) and Pic(${\mathbb F}$).
That is, a line bundle $\mathcal{X}$ on $X_0$ is equivalent to a line bundle 
$\mathcal{X}_{\mathbb P}$ on ${\mathbb P}$ and a line bundle $\mathcal{X}_{\mathbb F}$ on 
${\mathbb F}$ which agree when restricted to $R$.  
This means we must have that $\mathcal{X}_{\mathbb P} \cong 
\mathcal{O}_{\mathbb P}(d)$ and $\mathcal{X}_{\mathbb F} \cong 
\mathcal{O}_{\mathbb F}(cH-dE)$ for some $c$ and $d$.  Denote this bundle on 
$X_0$ by $\mathcal{X}(c, c-d)$.

We also have ${\mathbb P}$ and ${\mathbb F}$ as divisors on $X$, and the 
corresponding bundles $\mathcal{O}_X({\mathbb P})$ and 
$\mathcal{O}_X({\mathbb F})$.  The bundle $\mathcal{O}_X({\mathbb P})$ is 
disjoint from the fibers $X_t$ for $t\neq0$, but restricts to ${\mathbb P}$ 
as $\mathcal{O}_{\mathbb P}(-1)$ and restricts to ${\mathbb F}$ as 
$\mathcal{O}_{\mathbb F}(E)$.

Denote by $\mathcal{O}_X(d)$ the line bundle 
$\pi_2^*(\mathcal{O}_{{\mathbb P}^2}(d))$.  For $t\neq0$ the restriction of 
$\mathcal{O}_X(d)$ to $X_t$ is isomorphic to $\mathcal{O}_{{\mathbb P}^2}(d)$.
The restriction of $\mathcal{O}_X(d)$ to $X_0$ is the bundle 
$\mathcal{X}(d, 0)$.  The bundle $\mathcal{X}(d, 0)$ restricts to 
${\mathbb P}$ as $\mathcal{O}_{\mathbb P}(d)$ and to ${\mathbb F}$ as 
$\mathcal{O}_{\mathbb F}(dH-dE)$.

Let $\mathcal{O}_X(d, k)$ be the line bundle $\mathcal{O}_X(d) \otimes 
\mathcal{O}_X(k{\mathbb P})$.  When $\mathcal{O}_X(d, k)$ is restricted to $X_t$ for $t\neq0$
the result is isomorphic to $\mathcal{O}_{{\mathbb P}^2}(d)$, as before, but the restriction to
$X_0$ has changed.  The restriction of $\mathcal{O}_X(d, k)$ to $X_0$ is isomorphic to 
$\mathcal{X}(d, k)$, which restricts to ${\mathbb P}$ as $\mathcal{O}_{\mathbb P}(d-k)$ and 
restricts to ${\mathbb F}$ as $\mathcal{O}_{\mathbb F}(dH - (d-k)E)$.

In this way, all of the bundles $\mathcal{X}(d, k)$ on $X_0$ are seen as flat limits of 
the bundles $\mathcal{O}_{{\mathbb P}^2}(d)$ on the general fiber $X_t$ of this degeneration.
This is part of the flexibility this degeneration provides.  The rest of the flexibility lies
in the position of the points.

Take integers $n$ and $b$ such that $0 \leq b \leq n$, consider $n-b+1$ general points $p_0, 
p_1, \dots ,p_{n-b}$ in ${\mathbb P}$ and $b$ general points $p_{n-b+1} \dots ,p_n$ in 
${\mathbb F}$.  These can be realized as the limits of $n+1$ general points 
$p_{0,t}, p_{1,t} \dots,p_{n,t}$ in $X_t$.  Then we have the linear systems 
$\mathcal{L}_t (d, m_0, n, m) = \mathcal{L}(d, m_0, n, m)$ in
$X_t \cong {\mathbb P}_2$ for $t \neq 0$.  

On $X_0$ we have the system formed by sections of $\mathcal{X}(d, k)$ having a point of 
multiplicity $m_0$ at $p_0$ and multiplicity $m$ at points $p_1, \dots, p_n$.  This system 
will be called $\mathcal{L}_0 := \mathcal{L}_0 (d, k, m_0, n, b, m)$.  Any one of these systems
may be considered as the flat limit on $X_0$ of the system $\mathcal{L}_t = \mathcal{L}(d, m_0,
n, m)$.  We will say that $\mathcal{L}_0$ is obtained from $\mathcal{L}$ by a 
$(k, b)$ {\it degeneration}.

The linear system $\mathcal{L}_0$ restricts to ${\mathbb P}$ as the system 
$\mathcal{L}_{\mathbb P} := \mathcal{L}(d-k, m_0, n-b, m)$ and restricts to ${\mathbb F}$ as
the system $\mathcal{L}_{\mathbb F} := \mathcal{L}(d, d-k, b, m)$.
Divisors in the linear system $\mathcal{L}_0$ come in three types.  The first type consists
of a divisor $C_{\mathbb F}$ on ${\mathbb F}$ in the system $|dH - (d-k)E|$ and a divisor
$C_{\mathbb P}$ on ${\mathbb P}$ in the system $|(d-k)H|$, both of which satisfy the multiple 
point conditions, and which restrict to the same divisor on the curve $R$.

The second type is a divisor corresponding to a section of the bundle which is identically
zero on ${\mathbb P}$, and gives a divisor in the system $\mathcal{L}_{\mathbb F}$ which
contains the exceptional curve $E$ as a component.  A divisor in $\mathcal{L}_{\mathbb F}$
which contains $E$ is an element of the system $E + \mathcal{L}(d, d-k+1, b, m)$.  Since we are 
interested only in the dimension of this kernel system, we will denote it by 
$\hat{\mathcal{L}}_{\mathbb F}$ and (abusing notation) write $\hat{\mathcal{L}}_{\mathbb F} = 
\mathcal{L}(d, d-k+1, b, m)$.

The third type is similar.  It corresponds to a section of the bundle which is identically
zero on ${\mathbb F}$, and gives a divisor in the system $\mathcal{L}_{\mathbb P}$ which
contains the line $L$ as a component.  That is, it comes from an element of the system
$L + \mathcal{L}(d-k-1, m_0, n-b, m)$.  We will denote this kernel system by 
$\hat{\mathcal{L}}_{\mathbb P}$ and abuse notation further to write 
$\hat{\mathcal{L}}_{\mathbb P} = \mathcal{L}(d-k-1, m_0, b, m)$.

The four main linear systems are collected in the following table.  They are all quasi 
homogeneous, and in our case have $m=4$.  They also have smaller data than the original 
system, and give us the chance to argue by induction. 
$$\vbox{

\tabskip = 15pt
\halign{
        \hfil#&
        \hfil#&
        \hfil#&
        \hfil#&
        #\hfil

\cr
        $\hat{\mathcal{L}}_{\mathbb F}$ : & $d$ & $d-k+1$ & $b$ & 4 \cr
        $\mathcal{L}_{\mathbb F}$ : & $d$ & $d-k$ & $b$ & 4 \cr
        $\mathcal{L}_{\mathbb P}$ : & $d-k$ & $m_0$ & $n-b$ & 4 \cr
        $\hat{\mathcal{L}}_{\mathbb P}$ : & $d-k-1$ & $m_0$ & $n-b$ & 4 \cr  }
   
}$$

The following notation will be used.
$$\vbox{

\tabskip = 15pt
\halign{
        #\hfil&
        #\hfil

\cr
        $v$ & the virtual dimension of the general system, \cr
        $v_{\mathbb P}$ & the virtual dimension of the system on ${\mathbb P}$, \cr
        $v_{\mathbb F}$ & the virtual dimension of the system on ${\mathbb F}$, \cr
        $\hat{v}_{\mathbb P}$ & the virtual dimension of the kernel system on ${\mathbb P}$, \cr 
        $\hat{v}_{\mathbb F}$ & the virtual dimension of the kernel system on ${\mathbb F}$, \cr
        $\ell$ & the dimension of the general system, \cr
        $\ell_{\mathbb P}$ & the dimension of the system on ${\mathbb P}$, \cr
        $\ell_{\mathbb F}$ & the dimension of the system on ${\mathbb F}$, \cr 
        $\hat{\ell}_{\mathbb P}$ & the dimension of the kernel system on ${\mathbb P}$,  \cr 
        $\hat{\ell}_{\mathbb F}$ & the dimension of the kernel system on ${\mathbb F}$, and \cr
        $\ell_0$ & the dimension of the system $\mathcal{L}_0$. \cr    }

}$$
We have that $\ell \leq \ell_0$ by semi-continuity, and we will attempt to exploit the
inequality $v \leq \ell \leq \ell_0$ to show that the linear system $\mathcal{L}$ has the 
expected dimension.

Note the following relationships between the virtual dimensions of the main linear systems.
\begin{lemma} \label{ids} The following are identities in variables $d$, $m_0$, $n$, 
$k$, and $b$. \\ 
\indent a. $v_{\mathbb P} + v_{\mathbb F} = v+d-k.$ \\
\indent b.  $\hat{v}_{\mathbb P} + v_{\mathbb F} = v-1$. \\
\indent c.  $v_{\mathbb P} + \hat{v}_{\mathbb F} = v-1$. \\
\end{lemma}

\section{Classification of Special Systems with $m=4$}

One of the main results of \cite{miranda} is the following computation of the dimension of
$\mathcal{L}_0$.
\begin{theorem}  \label{Larry} Let $r_{\mathbb P} = \ell_{\mathbb P} -
\hat{\ell}_{\mathbb P} - 1$ and $r_{\mathbb F} = \ell_{\mathbb F} -
\hat{\ell}_{\mathbb F} - 1$.  Then \\
(a) If $r_{\mathbb P} + r_{\mathbb F} \leq d-k-1$, then
$\ell_0 = \hat{\ell}_{\mathbb P} + \hat{\ell}_{\mathbb F} + 1.$
(b) If $r_{\mathbb P} + r_{\mathbb F} \geq d-k-1$, then
$\ell_0 = \ell_{\mathbb P} + \ell_{\mathbb F} -d+k.$
\end{theorem}
The dimension computed in part (b) of the Theorem is the virtual dimension of the system
$\mathcal{L}$ by Lemma \ref{ids}, and will be used to show that a non-empty system has the 
expected dimension.  Part (a) of the Theorem is more useful for showing that a system is
empty.  The next two lemmas make this explicit.

\begin{lemma} \label{neglemma} Let $\mathcal{L}$ be a quasi homogeneous 
linear system with negative virtual dimension, $v \leq -1$.
If integers $k$ and $b$ can be found such that when a $(k, b)$
degeneration is executed \\ \indent (a) The systems 
$\mathcal{L}_{\mathbb F}$ and 
$\mathcal{L}_{\mathbb P}$ are both
non-special, and \\ \indent (b) the kernel systems 
$\hat{\mathcal{L}}_{\mathbb F}$ and $\hat{\mathcal{L}}_{\mathbb P}$ are
empty, \\ then $\mathcal{L}$ is empty.  \end{lemma}

\begin{proof}  If either $\mathcal{L}_{\mathbb F}$ or 
$\mathcal{L}_{\mathbb P}$ is empty then 
$\mathcal{L}$ is empty as well since the kernel systems are empty.  If 
$\mathcal{L}_{\mathbb F}$ and $\mathcal{L}_{\mathbb P}$ are not empty, then 
$$r_{\mathbb P} + r_{\mathbb F} = \ell_{\mathbb P} + 
\ell_{\mathbb F} = v_{\mathbb P} + v_{\mathbb F} = v+d-k \leq d-k-1.$$ 
The first equality follows from (b).  The
second is true because the systems are non-special and not empty.  The third
 equality is Lemma \ref{ids} a.
The final inequality holds by assumption $v \leq -1$.  Therefore
Theorem \ref{Larry} (a) applies and $\ell_0 = \hat{\ell}_{\mathbb P} +
\hat{\ell}_{\mathbb F} + 1 = -1$.  Now $\mathcal{L}$
must be empty since
$ -1 \leq \ell \leq \ell_0 = -1.$
\end{proof}

\begin{lemma} \label{poslemma} Let $\mathcal{L}$ be a quasi homogeneous 
linear system with  virtual dimension $v \geq -1$.
If integers $k <d$ and $b$ can be found such that when a $(k, b)$
degeneration is executed \\ \indent (a) The systems 
$\hat{\mathcal{L}}_{\mathbb F}$, $\mathcal{L}_{\mathbb F}$, 
$\mathcal{L}_{\mathbb P}$, and $\hat{\mathcal{L}}_{\mathbb P}$ are all
non-special, and \\ \indent (b) the systems 
$\mathcal{L}_{\mathbb F}$ and $\mathcal{L}_{\mathbb P}$ have virtual dimension at least -1,
 \\ then $\mathcal{L}$ has the expected dimension.  \end{lemma}

\begin{proof}  The proof relies on the identities from Lemma \ref{ids}. 
We claim that with the given hypotheses, $\hat{\ell}_{\mathbb P} +
\hat{\ell}_{\mathbb F} \leq v-1$.  There are three possibilities.
If both $\hat{\mathcal{L}}_{\mathbb P}$ and $\hat{\mathcal{L}}_{\mathbb
F}$ are empty, then $\hat{\ell}_{\mathbb P} +
\hat{\ell}_{\mathbb F} = -2 \leq v-1$ since $v\geq -1$.  If both systems
are non-empty and non special, then $\hat{\ell}_{\mathbb P} = 
\hat{v}_{\mathbb P}$ and $\hat{\ell}_{\mathbb F} = \hat{v}_{\mathbb F}$.
Then using the three identities we get 
$$\hat{\ell}_{\mathbb P} + \hat{\ell}_{\mathbb F} = \hat{v}_{\mathbb P}
+ \hat{v}_{\mathbb F} = v - 1 - d + k \leq v-1.$$  The inequality holds
since $k < d$.  If one of the systems is empty and the other is not,
identities b) and c) give  
$$\hat{\ell}_{\mathbb P} + \hat{\ell}_{\mathbb F} = -1 + \hat{v}_{\mathbb F} 
= v-1 -1-v_{\mathbb P} \leq v-1,$$
or
$$\hat{\ell}_{\mathbb P} + \hat{\ell}_{\mathbb F} = \hat{v}_{\mathbb P} -1 
= v-1 -1-v_{\mathbb F} \leq v-1.$$  The inequalities follow from
hypothesis (b).  Now,  
$$\vbox{
\tabskip = 15pt
\halign{
        \hfil#&
        \hfil#&
        #\hfil&
        #\hfil

\cr
        $r_{\mathbb P} + r_{\mathbb F}$ & = &  $\ell_{\mathbb P} - \hat{\ell}_{\mathbb
P} - 1 + \ell_{\mathbb F} - \hat{\ell}_{\mathbb F} -1$ &  \cr
         & = & $v_{\mathbb P} + v_{\mathbb F} - \hat{\ell}_{\mathbb P} 
- \hat{\ell}_{\mathbb F} - 2$ & by hypothesis \cr  
         & = & $ v+d-k- (\hat{\ell}_{\mathbb P} + \hat{\ell}_{\mathbb F})-2$ & 
by Lemma \ref{ids} a) \cr 
         & $\geq$ & $v+d-k-(v-1)-2$ & by the claim \cr
         & = & $d-k-1.$ & \cr }
   
}$$
We apply Theorem \ref{Larry} (b) to get $\ell_0 = \ell_{\mathbb P} +
\ell_{\mathbb F} -d +k = v_{\mathbb P} + v_{\mathbb F} -d+k = v$, by
Lemma \ref{ids} a.  Finally, we have $v \leq \ell \leq \ell_0 = v$.  Therefore
$\ell = v$ and $\mathcal{L}$ is non special.
\end{proof}

These are the basic tools in the proof of the main theorem.
\begin{theorem} A system $\mathcal{L}(d, m_0, n, 4)$ is special if and only 
if it is a (-1) special system, i.e., it is one of the systems listed in 
Theorem \ref{mosp}. \end{theorem}
\begin{proof} We may assume that $m_0 \leq d-6$ and that $d\geq 6$.  
The proof is by induction on $d$.  Assume the theorem is true for smaller 
values of $d$.  Then assume that $\mathcal{L}$ is not (-1) special, and 
prove that it is non-special.  

Begin with the case $v \leq -1$.  We perform a $(k, b)$ degeneration with 
$k=3$.  This gives the following relevant systems. 
$$\vbox{

\tabskip = 15pt
\halign{
        \hfil#&
        \hfil#&
        \hfil#&
        \hfil#&
        #\hfil

\cr
        $\hat{\mathcal{L}}_{\mathbb F}$ : & $d$ & $d-2$ & $b$ & 4 \cr
        $\mathcal{L}_{\mathbb F}$ : & $d$ & $d-3$ & $b$ & 4 \cr
        $\mathcal{L}_{\mathbb P}$ : & $d-3$ & $m_0$ & $n-b$ & 4 \cr
        $\hat{\mathcal{L}}_{\mathbb P}$ : & $d-4$ & $m_0$ & $n-b$ & 4 \cr  }
   
}$$
We wish to find a $b$ so that both $\hat{\mathcal{L}}_{\mathbb P}$  and 
$\hat{\mathcal{L}}_{\mathbb F}$ are empty, and all four systems are 
non-special.  We pick $b > \frac{d}{3}$ so that 
$\hat{\mathcal{L}}_{\mathbb F}$ is not (-1) special, and thus non-special 
by section \ref{largem_0}.  Then $b > \frac{3d-1}{10}$, which makes 
$\hat{v}_{\mathbb F} \leq -1$.  Therefore $\hat{\mathcal{L}}_{\mathbb F}$ 
is empty.  If, in addition, $b \leq \frac{4d-2}{10}$, then 
$\hat{v}_{\mathbb P} \leq -1$ since $v - \hat{v}_{\mathbb P} = 4d - 2 +10b$ 
and $v \leq -1$.  Now if $\hat{\mathcal{L}}_{\mathbb P}$ is not (-1) special 
(and non-special by induction), it will be empty.  If $m_0 < d-6$ this can 
be ensured by choosing $b$ so that $n-b$ is odd.  This is also enough to 
conclude that $\mathcal{L}_{\mathbb P}$ is not (-1) special (and 
non-special by induction).  In order for $\mathcal{L}_{\mathbb F}$ to be 
not (-1) special (and thus non-special by section \ref{largem_0}), we choose 
$b < \frac{4d-4}{10}$.

An integer $b$, satisfying the inequalities
$$\frac{3d-1}{10} < \frac{d}{3} < b < \frac{4d-4}{10} < \frac{4d-2}{10},$$
and such that $n-b$ is odd, may be found when $d$ is 29, 32, 34, 35, or 
$d \geq 37$.   This choice of $b$ makes all of the systems involved 
non-special and the kernel systems empty.  We apply Lemma
\ref{neglemma}.
This proves that $\mathcal{L}$ is
empty, provided $d \geq 37$, $m_0 < d-6$ and the conjecture holds for 
lower $d$.

If $m_0 = d-6$, $\hat{\mathcal{L}}_{\mathbb P}$ might be (-1) special even 
when $n-b$ is odd,
but only if $n-b \leq \frac{3(d-4)-1}{9}$.  If 
$\hat{\mathcal{L}}_{\mathbb P}$ is not (-1) special, Lemma \ref{neglemma} can 
be used to conclude that $\mathcal{L}$ is empty.
We would like $$n > \frac{3(d-4)-1}{9} + b.$$
For a given $d$ and $m_0$, we only need to prove the theorem for the
smallest $n$ which makes $v$ negative.  For $m_0 = d-6$, this value is
the smallest integer $n$ such that $n > \frac{7d-15}{10}$.  It is clear
that we need to choose $b$ as small as possible.  We will pick
$\frac{d}{3} < b \leq \frac{d}{3} +2$ such that $n-b$ is odd.  
Then $\frac{3(d-4)-1}{9} + b
\leq \frac{6d+5}{9}$.  Now $$n > \frac{7d-15}{10} > \frac{6d+5}{9}$$
whenever $d \geq 62$.  This guarantees that $\hat{\mathcal{L}}_{\mathbb
P}$ is not (-1) special when $m_0 = d-6$ and $d\geq62$.  For values of
$d$ between 62 and 39, one may directly check that this choice of $b$ makes 
$\hat{\mathcal{L}}_{\mathbb P}$ not (-1) special for the smallest value of
$n$.  The same is true for $d = 38, 37, 35,$ and 34.  We are assuming the 
theorem for smaller values of $d$, so Lemma 
\ref{neglemma} applies and these systems are empty.  The only
outstanding case with $d\geq 37$ is $d=39$ and $m_0 = 33$.

For $d=39$ and $m_0 = 33$, the smallest $n$ which makes $v$ negative is 
$n=26$.  We choose $k=3$ and $b=13$.  This makes 
$\hat{\mathcal{L}}_{\mathbb F}$ non-special and empty by section 
\ref{largem_0} and $\hat{\mathcal{L}}_{\mathbb P}$ non-special and empty
 by induction.  $\mathcal{L}_{\mathbb F}$ is not (-1) special since $b <
\frac{4d-4}{10}$ and non-special by section \ref{largem_0}.  
$\mathcal{L}_{\mathbb P}$ is not (-1) special since $n-b$ is odd for
$n=26$ and non-special by induction.
Now all the systems are non-special and the kernel systems are empty so we use 
Lemma \ref{neglemma} to conclude that this system is empty.
  We have proven the theorem
for all $d \geq 37$ provided it is true for smaller values of $d$.  The
theorem will be proved case by case for smaller values of $d$.

If $d=36$, we perform a $(3, 13)$ degeneration giving the systems
$$\vbox{

\tabskip = 15pt
\halign{
        \hfil#&
        \hfil#&
        \hfil#&
        \hfil#&
        #\hfil

\cr
        $\hat{\mathcal{L}}_{\mathbb F}$ : & 36 & 34 & 13 & 4 \cr
        $\mathcal{L}_{\mathbb F}$ : & 36 & 33 & 13 & 4 \cr
        $\mathcal{L}_{\mathbb P}$ : & 33 & $m_0$ & $n-13$ & 4 \cr
        $\hat{\mathcal{L}}_{\mathbb P}$ : & 32 & $m_0$ & $n-13$ & 4 \cr  }
   
}$$
We are assuming that $m_0 \leq 30$.
If all these systems are non-special and the kernel systems are empty
then we use Lemma \ref{neglemma} to conclude that $\mathcal{L}$ is empty.  
The kernel systems have negative
virtual dimension, so are empty if they are non-special.
The systems $\hat{\mathcal{L}}_{\mathbb F}$ and $\mathcal{L}_{\mathbb
F}$ are not special by section \ref{largem_0}.  
The system $\hat{\mathcal{L}}_{\mathbb P}$ is (-1)
special only if it has the form $\mathcal{L}(32, 28, 16, 4)$ or 
$\mathcal{L}(32, 30, x, 4)$ with $x \leq 10$.  The first comes from the 
system 
$\mathcal{L} = \mathcal{L}(36, 28, 29, 4)$ which has virtual dimension
$v=6$, contrary to hypothesis.  The second type comes from a system 
$\mathcal{L} = \mathcal{L}(36, 30, x+13, 4)$ with virtual dimension
$v=107-10x$.  Again, $v$ is positive (since $x \leq 10$), contrary to 
hypothesis.  The
system $\mathcal{L}_{\mathbb P}$ is (-1) special only if it has the form
$\mathcal{L}(33, 29, 16, 4)$.  This comes from the system $\mathcal{L} =
\mathcal{L}(36, 29, 29, 4)$ with virtual dimension $v=-23$.  The system 
$\mathcal{L}_{\mathbb P} = \mathcal{L}(33, 29, 16, 4)$ is special, but 
Cremona reduces to 
the class of a line, and so has dimension 2.  $\mathcal{L}_{\mathbb F}$
is non-special with virtual dimension 11.
So in this case, $r_{\mathbb P} + r_{\mathbb F} = \ell_{\mathbb P} + 
\ell_{\mathbb F} = 2 + 11 < d-k-1 = 32$ and we may appeal directly to 
Theorem \ref{Larry} (a) to conclude that $\mathcal{L}$ is empty.  

The next open case is $d=33$. (The cases $d=34$ and $d=35$ were
mentioned earlier.)  We perform a $(3, 11)$ degeneration giving the systems
$$\vbox{

\tabskip = 15pt
\halign{
        \hfil#&
        \hfil#&
        \hfil#&
        \hfil#&
        #\hfil

\cr
        $\hat{\mathcal{L}}_{\mathbb F}$ : & 33 & 31 & 11 & 4 \cr
        $\mathcal{L}_{\mathbb F}$ : & 33 & 30 & 11 & 4 \cr
        $\mathcal{L}_{\mathbb P}$ : & 30 & $m_0$ & $n-11$ & 4 \cr
        $\hat{\mathcal{L}}_{\mathbb P}$ : & 29 & $m_0$ & $n-11$ & 4 \cr  }
   
}$$
The kernel systems have negative virtual dimension, hence they are empty if 
they are non-special.  $\hat{\mathcal{L}}_{\mathbb F}$ and 
$\mathcal{L}_{\mathbb F}$ are non-special by section \ref{largem_0}.  If 
$\mathcal{L}_{\mathbb P}$ and $\hat{\mathcal{L}}_{\mathbb P}$ are non-
special we apply Lemma \ref{neglemma} to conclude that $\mathcal{L}$ is
empty.  $\hat{\mathcal{L}}_{\mathbb P}$ is only (-1) special if it of
the form $\mathcal{L}(29, 25, 14, 4)$ or $\mathcal{L}(29, 27, x, 4)$
where $x \leq 9$.  Both of these come from systems with positive virtual
dimension.  $\mathcal{L}_{\mathbb P}$ is (-1) special only if it is
$\mathcal{L}(30, 27, 12, 4)$ or $\mathcal{L}(30, 26, 14, 4)$.  In these
cases we appeal directly to Theorem \ref{Larry} (a).

In both cases $\mathcal{L}_{\mathbb F}$ is non-special of dimension 19.
$\mathcal{L}(30, 27, 12, 4)$ Cremona reduces to the zero dimensional
space of constant polynomials.  Now $r_{\mathbb P} + r_{\mathbb F} = 
\ell_{\mathbb P} + \ell_{\mathbb F} = 0 + 19 < d-k-1 = 29$, so we use Theorem 
\ref{Larry} (a) to conclude that $\mathcal{L}$ is empty.  
In the other case, we notice that $\mathcal{L}(30, 26, 14, 4)$ Cremona 
reduces to the 5 dimensional space of quadratics.  This means that 
$r_{\mathbb P} + r_{\mathbb F} = 
\ell_{\mathbb P} + \ell_{\mathbb F} = 5 + 19 < d-k-1 = 29$, and Theorem 
\ref{Larry} (a) tells us that $\mathcal{L}$ is empty.

When $d=32$ we pick $b$ as before ($\frac{d}{3} < b < \frac{4d-4}{10}$
and $n-b$ odd).  This allows us to apply Lemma \ref{neglemma} 
unless $\mathcal{L} = \mathcal{L}(32, 26, n, 4)$.  It is enough to 
prove $\mathcal{L}(32, 26, 21, 4)$ is empty.  We perform a $(3, 11)$ 
degeneration in this case and apply Lemma
\ref{neglemma}.

For $d=31$, a $(3, 11)$ degeneration produces systems which are not (-1) 
special and kernel systems which are empty, unless $\mathcal{L} =
\mathcal{L}(31, 24, 25, 4)$ or $\mathcal{L}(31, 25, x+11, 4)$ with
$x\leq 8$.  The later does not have negative virtual dimension.  The
former degenerates as
$$\vbox{

\tabskip = 15pt
\halign{
        \hfil#&
        \hfil#&
        \hfil#&
        \hfil#&
        #\hfil

\cr
        $\hat{\mathcal{L}}_{\mathbb F}$ : & 31 & 29 & 11 & 4 \cr
        $\mathcal{L}_{\mathbb F}$ : & 31 & 28 & 11 & 4 \cr
        $\mathcal{L}_{\mathbb P}$ : & 28 & 24 & 14 & 4 \cr
        $\hat{\mathcal{L}}_{\mathbb P}$ : & 27 & 24 & 14 & 4 \cr  }
   
}$$
$\mathcal{L}_{\mathbb F}$ is non-special of dimension $\ell_{\mathbb F}
= 11$.  $\mathcal{L}_{\mathbb P}$ Cremona reduces to the zero dimensional 
class of a quadruple line.  In this case, $r_{\mathbb P} + r_{\mathbb F} = 
\ell_{\mathbb P} + \ell_{\mathbb F} = 11+0 < d-k-1 = 27$, and Theorem 
\ref{Larry} (a) tells us that $\mathcal{L}$ is empty.

If $d=30$, we perform a $(3, 11)$ degeneration.  This produces systems
which are not (-1) special and kernel systems which are empty unless
$\hat{\mathcal{L}}_{\mathbb P} = \mathcal{L}(26, 22, 12, 4),$ 
$\mathcal{L}(26, 23, 10, 4),$ or $\mathcal{L}(26, 24, x, 4)$, where
$x\leq 8$.  These exceptions all come from systems with positive virtual
dimension, therefore Lemma \ref{neglemma} handles all cases with $d=30$.

When $d=29$ we pick $b$ as before ($\frac{d}{3} < b < \frac{4d-4}{10}$
and $n-b$ odd).  This satisfies the conditions of Lemma \ref{neglemma}
unless $\mathcal{L} = \mathcal{L}(29, 23, n, 4)$.  It is enough to show
that $\mathcal{L}(29, 23, 19, 4)$ is empty.  A $(3, 10)$ degeneration
allows us to use Lemma \ref{neglemma} in this instance as well.

For $d=28$, a $(3, 10)$ degeneration satisfies the conditions of Lemma 
\ref{neglemma} in all but five cases, only two of which have negative
virtual dimension.  In these cases we use Cremona reduction to find the
dimension of $\mathcal{L}_{\mathbb P}$ and apply Theorem \ref{Larry} (a).
The $d=27$ case has one exception when a $(3, 9)$ degeneration is used,
and it is handled in the same way.  A $(3, 9)$ degeneration also suffices to
prove the theorem in case $d=26$ and 25.

When $d=24$, a $(3, 9)$ degeneration works unless $\mathcal{L}$ is 
$\mathcal{L}(24, 18, 17, 4)$, $\mathcal{L}(24, 17, 19, 4)$, or 
$\mathcal{L}(24, 16, 19, 4)$.  In the first two of these exceptions
$\mathcal{L}_{\mathbb P}$ is special and we proceed as above using
Cremona reduction to find the dimension.  In the final case, it is 
$\hat{\mathcal{L}}_{\mathbb P}$ which is special and we must find
another approach.  A $(4, 10)$ degeneration allows us to
use Lemma \ref{neglemma} to conclude that $\mathcal{L}$ is empty in this
case as well.

If $d=23$, we use a $(3, 8)$ degeneration.  This lets us apply Lemma
\ref{neglemma} in all but two cases.  In these cases, we use Cremona
reduction and Theorem \ref{Larry} (a).  For $d=22$, a $(3, 8)$ degeneration
allows us to apply Lemma \ref{neglemma} unless $\mathcal{L} = \mathcal{L}(22,
16, 14, 4)$.  This system Cremona reduces to the system $\mathcal{L}(8,
0, 15, 2)$.  The theorem is true for $m=2$ (\cite{miranda}), and this
system is not (-1) special, so this system is empty.  Therefore
$\mathcal{L}$ is empty, as well.

When $18 \leq d \leq 21$, we use a $(3, 7)$ degeneration.  We may appeal
to Lemma \ref{neglemma} in all but a finite number of cases.  In all but one of
these cases we use Cremona reduction to determine the dimension of 
$\mathcal{L}_{\mathbb P}$ and apply Theorem \ref{Larry} (a).  In case
$\mathcal{L} = \mathcal{L}(19, 12, 15, 4)$, we use a $(3, 8)$
degeneration and Lemma \ref{neglemma}.  If $d=17$, a $(3, 6)$
degeneration succeeds.  If $d=16$ there are eleven possibilities for $m_0$.
  In each case, it is enough to prove the theorem for the smallest $n$
 which makes $v$ negative and the system not (-1) special.  We have the 
following systems

$$\vbox{

\tabskip = 15pt
\halign{
        \hfil#&
        \hfil#&
        \hfil#&
        #\hfil

\cr
        $\mathcal{L}(16,$ & 10, & 10, & 4)  \cr
        $\mathcal{L}(16,$ & 9, & 11, & 4)   \cr
        $\mathcal{L}(16,$ & 8, & 12, & 4) \cr
        $\mathcal{L}(16,$ & 7, & 13, & 4) \cr  
        $\mathcal{L}(16,$ & 6, & 14, & 4) \cr  
        $\mathcal{L}(16,$ & 5, & 14, & 4)   \cr
        $\mathcal{L}(16,$ & 4, & 15, & 4)  \cr  
        $\mathcal{L}(16,$ & 3, & 15, & 4)  \cr
        $\mathcal{L}(16,$ & 2, & 15, & 4)  \cr
        $\mathcal{L}(16,$ & 1, & 16, & 4)  \cr
        $\mathcal{L}(16,$ & 0, & 16, & 4)  \cr }
   
}$$

\noindent The first system Cremona reduces to a homogeneous system with $m=2$, which is 
empty \cite{miranda}.  A $(3, 6)$ degeneration is used for all but the third
 and ninth systems.  For these degenerations, $\mathcal{L}_{\mathbb F}$ is 
special but we use Cremona transformations to find its dimension and apply 
Theorem \ref{Larry} (a).  The remaining systems yield to a $(4, 7)$ 
degeneration and Lemma \ref{neglemma}.

For $d=15$ or 14, a $(3, 5)$
degeneration works.  A $(3, 5)$ degeneration also works when $d=13$ in
all but one case.  The case $\mathcal{L} = \mathcal{L}(13, 5, 9, 4)$
did not yield to any of the approaches used so far.  To see that this system is 
empty, Maple was used to analyze the $105\times105$ matrix constructed using the 
points $p_0 = (0, -3)$, $p_1 = (8, 3)$, $p_2 = (4, -4)$, $p_3 = (-5, -5)$,
$p_4 = (-5, -2)$, $p_5 = (3, -1)$, $p_6 = (-5, -9)$, $p_7 = (8, 5)$, 
$p_8 = (5, 8)$, and $p_9 = (-1, 4)$.  The matrix was found to have full rank,
therefore there are no thirteenth degree polynomials passing through these
points with multiplicities.  This implies that $\mathcal{L}(13, 5, 9, 4)$ is 
empty for points in general position.

In the $d=12$ case, we again list the systems corresponding to possible 
values of $m_0$ and the critical $n$ for each.  We have
the following seven systems.

$$\vbox{

\tabskip = 15pt
\halign{
        \hfil#&
        \hfil#&
        \hfil#&
        #\hfil

\cr
        $\mathcal{L}(12,$ & 6, & 8, & 4)  \cr
        $\mathcal{L}(12,$ & 5, & 8, & 4)   \cr
        $\mathcal{L}(12,$ & 4, & 9, & 4) \cr
        $\mathcal{L}(12,$ & 3, & 9, & 4) \cr  
        $\mathcal{L}(12,$ & 2, & 9, & 4) \cr  
        $\mathcal{L}(12,$ & 1, & 9, & 4)   \cr
        $\mathcal{L}(12,$ & 0, & 10, & 4)  \cr   }
   
}$$

\noindent It is sufficient to show that the second, sixth, and seventh systems are
empty, as these imply the rest are empty.  The second system Cremona
reduces to $\mathcal{L}(8, 1, 8, 3)$, which is non-special and empty by
\cite{miranda}.  To show the sixth is empty, we employ a $(4, 5)$
degeneration and Lemma \ref{neglemma}.  Finally, we use a $(4, 5)$
degeneration with the seventh system.  $\mathcal{L}_{\mathbb P}$ Cremona
reduces to a quadruple line, and so has dimension 0.  We apply Theorem
\ref{Larry} (a).

The $d=11$ case proceeds similarly.  There are six possibilities $0
\leq m_0 \leq 5$, and we only need to check the smallest $n$ in each.
We only need to show that $\mathcal{L}(11, 4, 7, 4)$ and
$\mathcal{L}(11, 0, 8, 4)$ are empty since these imply the other four
are empty.  In both of these cases we use a $(3, 4)$ degeneration.
$\mathcal{L}_{\mathbb F}$ is special, but Cremona reduces to the two dimensional 
class of a line.  In both cases $\mathcal{L}_{\mathbb P}$ is not (-1)
special with virtual dimension 4.  We apply Theorem
\ref{Larry} (a).

In the $d=10$ case, for similar reasons, 
we only need to show that $\mathcal{L}(10, 3, 6, 4)$
and $\mathcal{L}(10, 0, 7, 4)$ are empty.  The first Cremona reduces to
the empty system of constant polynomials with a double point.  The second
reduces to the empty system of quadratics with a triple point.
When $d=9$, it suffices to show that $\mathcal{L}(9, 3, 5, 4)$ and
$\mathcal{L}(9, 0, 6, 4)$ are empty.  The former reduces to the empty
system of lines through three general points.  The latter reduces to the
empty system of non-zero
constant polynomials passing through three points.
For $d=8$, we need to show that $\mathcal{L}(8, 1, 5, 4)$ and
$\mathcal{L}(8, 0, 6, 4)$ are empty.  The first reduces to the empty
system of non-zero
constant polynomials passing through a point.  The second we may
conclude is empty by applying a $(3, 3)$ degeneration and Lemma
\ref{neglemma}.  Finally, we must show that $\mathcal{L}(7, 0, 4, 4)$
and $\mathcal{L}(6, 0, 4, 4)$ are empty.  $\mathcal{L}(7, 0, 4, 4)$
reduces to the empty system of quadratics with a quadruple point.  
$\mathcal{L}(6, 0, 4, 4)$ reduces to the empty system of constant
polynomials with a quadruple point.

Now we consider the $v \geq -1$ case.  Again, we assume that $m_0 \leq d-
6$ and that $d \geq 6$.  We assume that $\mathcal{L}$
 is not (-1) special with $v \geq -1$ and prove that it is non special.  For 
each $d$ and $m_0$ it is enough to prove the theorem for the largest $n$ 
which makes the virtual dimension $v\geq 0$.  
If $\mathcal{L}(d, m_0, n, 4)$ is non empty and non special, 
then 
$\mathcal{L}(d, m_0, n^\prime, 4)$ will be non-empty and non-special for 
$n^\prime < n$.  
This is because the conditions imposed on curves of degree $d$ in 
$\mathcal{L}(d, m_0, n^\prime, 4)$ are a subset of the independent conditions
 for the system 
$\mathcal{L}(d, m_0, n, 4)$.  For $m_0 = d-6$, $v \geq 0$ when $n \leq 
\frac{7d-15}{10}$.  For lower $m_0$, the largest value of $n$ which makes 
$v\geq 0$ must be at least as big as in the $m_0 = d-6$ case.  Therefore, we 
will assume that $n >  \frac{7d-15}{10} - 1$.

We use a $(3, b)$ degeneration where $b$ is chosen to satisfy the 
hypotheses of Lemma \ref{poslemma}.  $\hat{v}_{\mathbb F}$ is not (-1) 
special when $b > \frac{d}{3}$.  $v_{\mathbb F}$ is not (-1) special if 
$b < \frac{2d-2}{5}$.  $v_{\mathbb P}$ is not (-1) special when $n-b$ is odd.
  $\hat{v}_{\mathbb P}$ is not (-1) special if $n-b$ is odd and $m_0 
\neq d-6$.  Now $\hat{v}_{\mathbb F}$ and $v_{\mathbb F}$ are non-special and
 $v_{\mathbb P}$ and $\hat{v}_{\mathbb P}$ are non special by the inductive 
hypothesis.  We need $b \leq \frac{2d-1}{5}$ to make $v_{\mathbb F} \geq -1$.
  To force $v_{\mathbb P} \geq -1$, the second identity in the proof of Lemma
 \ref{poslemma} gives us that it is enough to have that $v\geq-1$ and 
$\hat{v}_{\mathbb F} \leq -1$.  The former is true by 
hypothesis and the later is true when $b> \frac{3d}{10}$.  All of this may 
be achieved (when $m_0 < d-6$) by choosing $$\frac{d}{3} < b < 
\frac{2d-2}{5}$$ so that $n-b$ is odd.  As before, we can do this for $d=29$,
 32, 34, 35, or whenever $d\geq37$.  Recall that $n >  \frac{7d-25}{10}$, so 
that this choice makes $b<n$ for all $d \geq 8$.

When $m_0 = d-6$, $\hat{v}_{\mathbb P}$ will be (-1) special if $n-b \leq 
\frac{3(d-4)-1}{9}$.  
 For $m_0 = d-6$ the largest $n$ making $v \geq
 0$ is the integer 
$\frac{7d-25}{10} < n \leq \frac{7d-15}{10}$.  We need $b < n- 
\frac{3(d-4)-1}{9}$.  Putting 
these last two inequalities together, we see that we should choose $b < 
\frac{33d-95}{90}$.  
This guarantees that $\hat{v}_{\mathbb P}$ is non-special even when $m_0 = 
d-6$.

If we select $\frac{d}{3} < b < \frac{33d-95}{90}$ such that $n-b$ is odd, 
then we may appeal 
to Lemma \ref{poslemma} to conclude that $\mathcal{L}$ is non-special.  This 
may be done for 
$d\geq 91$.  For $37 \leq d \leq 90$, the hypotheses of Lemma
\ref{poslemma} are satisfied (for $m_0=d-6$ and the largest $n$ which makes 
$v\geq-1$) 
if we choose $\frac{d}{3} \leq b \leq \frac{d}{3} +1$, for all but
$d=46$ and $d=40$.  For the system $\mathcal{L}(46, 40, 30, 4)$, we
perform a $(3, 16)$ degeneration.  $\hat{\mathcal{L}}_{\mathbb F}$ is
empty, $\mathcal{L}_{\mathbb F}$ is non-special of dimension 21,
$\mathcal{L}_{\mathbb P}$ is non-special (by induction) of dimension 29, and
$\hat{\mathcal{L}}_{\mathbb P}$ is expected to be empty, but is (-1)
special.  $\hat{\mathcal{L}}_{\mathbb P} = \mathcal{L}(42, 40, 14, 4)$
Cremona reduces to the zero dimensional space of constants.  We compute
$r_{\mathbb P} + r_{\mathbb F} = (29 - 0 - 1) + (21 -(-1) -1) = 49 > d-
k-1 = 46-3-1.$  Therefore, $\ell_0 = 21+29 - 46 +3 = 7$ by Theorem 
\ref{Larry} (b),
 and this is also the virtual dimension. 
 The $d=40$ case follows in the same fashion.  Use a $(3, 14)$ degeneration.  
$\hat{\mathcal{L}}_{\mathbb P}$ is special, but Cremona reduces to the zero 
dimensional space 
of constants.  The hypothesis of Theorem \ref{Larry} (b) are satisfied
and $v = \ell_0$ since $v_{\mathbb P} = \ell_{\mathbb P}$ and
$v_{\mathbb F} = \ell_{\mathbb F}$.

The theorem is proved for $d \geq 37$, provided it is true for smaller values
 of $d$.  
We now prove the theorem for $d < 37$ case by case.  In each case we choose
$b$ so that $v_{\mathbb F}$ and $v_{\mathbb F}$ are at least -1, then check 
that all the systems are not (-1) special.  If they are all not (-1) special
we may apply Lemma \ref{poslemma}.  When $d=36$ we use
a $(3, 13)$ degeneration.  We analyzed this degeneration already in the $v\leq
-1$ case.  If the four main systems are all non-special and 
$\mathcal{L}_{\mathbb P}$
and $\mathcal{L}_{\mathbb F}$ have dimension at least -1 we apply
Theorem \ref{Larry} (b) to conclude that $\mathcal{L}$ has the expected
dimension.  Now, $b$ was chosen so that $\mathcal{L}_{\mathbb P}$
and $\mathcal{L}_{\mathbb F}$ have virtual dimension at least -1, so it
is enough that the systems are all non-special.  We saw earlier that
this only fails (when $v \geq -1$) if $\mathcal{L} = \mathcal{L}(36, 28,
29, 4)$ or if $\mathcal{L}$ is of the form $\mathcal{L}(36, 30, x+13,
4)$ where $x\leq 10$.  In the second case it is enough to prove the
theorem for the largest $n$ making $v\geq-1$, so we may assume
$\mathcal{L} = \mathcal{L}(36, 30, 23, 4)$.  In both cases
$\hat{\mathcal{L}}_{\mathbb P}$ is special, but we may use Cremona
transformations to find its dimension.  Then directly apply Theorem
\ref{Larry} (b).

For $d=35$, we may choose $b$ to be 12 or 13 to make $n-b$ odd (for the largest
 $n$ making $v$ at least 0) and use Lemma \ref{poslemma}.  If $d=34$, this $b$
 also works unless $\mathcal{L} = \mathcal{L}(34, 28, 22, 4)$.  In this case we
appeal to Theorem \ref{Larry} (b) (using Cremona reduction to find the 
dimension of $\hat{\mathcal{L}}_{\mathbb P}$).  When $d=33$ we use a $(3, 11)$
degeneration and Lemma \ref{poslemma} unless $\mathcal{L} = \mathcal{L}(33,
25, 25, 4)$ or is of the form $\mathcal{L}(33, 27, x+11, 4)$ where $x \leq9$. 
In the first case we appeal to Theorem \ref{Larry} (b).  In the second case,
note that we only care about the largest $n$ (the largest $x$).  This is the
system $\mathcal{L} = \mathcal{L}(33, 27, 20, 4)$, but $\mathcal{L}(33, 27,
21, 4)$ is non-empty and non special (by Lemma \ref{poslemma}) so $\mathcal{L}$
is as well.

When $d=32$, we may choose $b = 11$ or 12 so that $n-b$ is odd for the largest
$n$ making $v$ at least 0 and use Lemma \ref{poslemma}.  For $d=31$ a $(3, 11)$
degeneration allows us to apply Lemma \ref{poslemma}.  A $(3, 11)$ degeneration
also allows us to apply Lemma \ref{poslemma} in all but three cases.  For
these exceptions, we use Cremona reduction to find the dimension of the 
special systems and use Theorem \ref{Larry} (b).  If $d=29$, we let $b$ be 
either 10 or 11 so that $n-b$ is odd and use Lemma \ref{poslemma} unless 
$\mathcal{L} = \mathcal{L}(29, 23, 18, 4)$.  In this case we use a $(3, 10)$
degeneration and Theorem \ref{Larry} (b).

For every $d$ from 17 to 28 the process is the same.  For each $m_0 \leq d-6$ and the 
corresponding largest $n$, there is a $b$ (between $\frac{3d}{10}$ and $\frac{2d-1}{5}$ so that
$v_{\mathbb F}$ and $v_{\mathbb P}$ are at least -1) such that either a $(3, b)$ degeneration
satisfies the hypothesis of Lemma \ref{poslemma}, or so that we may use Cremona reduction
to find the dimension of the special systems and apply Theorem \ref{Larry} (b).  

When $d=16$, we may use a $(3, 5)$ degeneration and Theorem \ref{Larry} (b) in all but two 
cases.  In these cases ($\mathcal{L} = \mathcal{L}(16, 6, 13, 4)$ or 
$\mathcal{L}(16, 1, 15, 4)$) a $(4, 7)$ degeneration satisfies the requirements of Lemma
\ref{poslemma}.  For $d=15$ or 14, a $(3, 5)$ degeneration works.  

If $d=13$ there are eight cases to check.  For $m_0 = 7$ or 6 the systems (with the largest $n$) 
Cremona reduce to a known case and are non-special.  Either a $(3, 4)$ or a $(3, 5)$ 
degeneration works for the rest of the cases except $m_0 = 2$.  To prove the theorem for the
system $\mathcal{L}(13, 2, 10, 4)$ a $(6, 7)$ degeneration satisfies the hypotheses of Lemma
\ref{poslemma}.  When $d=12$ there are six cases to check.  The top two ($m_0 = 5$ or 6)
Cremona reduce to a known case and are non-special.  A $(3, 4)$ or a $(4, 5)$ degeneration
works for the rest.

For $6 \leq d \leq 11$, $0 \leq m_0 \leq d-6$, and the corresponding largest $n$ making $v$ at
least 0 and the system non-special, the linear systems all Cremona reduce to known cases and
are non-special.

Therefore, all linear systems $\mathcal{L}(d, m_0, n, m)$ which are not (-1) special are non-
special.  In other words the only special quasi-homogeneous linear systems with $m=4$ are the 
(-1) special systems listed in Theorem \ref{mosp}.

\end{proof}

\end{document}